\documentclass[11pt]{article}
\usepackage{amssymb,amsmath,latexsym}
\usepackage{eufrak}
\usepackage{amsthm}
\author{Eric Emtander\footnote{Department of Mathematics, 
Stockholm University, 106 91 Stockholm, erice@math.su.se}}
\title{Betti numbers of hypergraphs}
\date{}
\begin{document}
\maketitle
\newtheorem{thm}{Theorem}[section]
\newtheorem{lemma}[thm]{Lemma}
\newtheorem{prop}[thm]{Proposition}
\newtheorem{cor}[thm]{Corollary}
\let\sse=\subseteq \let\ssm=\smallsetminus\let\hra=\hookrightarrow
\let\car=\curvearrowright\let\ss=\subset\let\t=\textrm\let\mc=\mathcal
\abstract{In this paper we study some algebraic properties of hypergraphs, in particular 
their Betti numbers. We define some different types of complete hypergraphs, which to 
the best of our knowledge, are not previously considered in the literature. Also, in a 
natural way, we define a product on hypergraphs, which in a sense is dual to 
the join operation on simplicial complexes. For such product, we give a general formula 
for the Betti numbers, which specializes neatly in case of linear resolutions.}

\section{Introduction}

Let ${\mc X}$ be a finite set and ${\mc E}=\{E_1,...,E_s\}$ a finite collection of 
non empty subsets of ${\mc X}$. The pair ${\mc{ H=(X,E)}}$ is called a {\bf hypergraph}. 
The elements of ${\mc X}$ are called the {\bf vertices} and the elements of ${\mc E}$ 
are called the {\bf edges} of the hypergraph. If we want to specify what hypergraph we 
consider, we may write $\mc{X(H)}$ and $\mc{E(H)}$ for the vertices and edges 
respectively.

The hypergraphs that we will consider, can all be seen as natural generalizations of 
the ordinary complete graph $K_n$, on $n$ vertices. Our main tools are familiar concepts 
in combinatorial algebra, such as Hochster's formula, the Mayer-Vietoris sequence and 
K\"unneth's tensor formula.

A hypergraph is called {\bf simple} if: (1) $|E_i|\ge2$ for all 
$i=1,...,s$ and (2) $E_j\sse E_i$ implies $i=j$. If the cardinality of ${\mc X}$ is 
$n$ we often just use the set $[n]=\{1,2,...,n\}$ instead of $\mc{X}$.

We frequently identify a vertex $v_i$ of $\mc H$ with 
a variable $x_i$ of a polynomial ring $k[x_1,...,x_n]$ over some field $k$, or 
with its corresponding characteristic vector $v(v_i)=(0,...,0,1,0,...,0)$ in 
$\mathbb{N}^n$, consisting of only zeros except in the $i$'th position were there is a 1. 
Hence we choose to consider 0 to be a natural number. This also allows us to identify a 
subset $V$ of $[n]$ with its characteristic vector $v(V)=\sum_{i\in V}v(v_i)$. 
We use bold letters to denote vectors and if 
${\bf w}=(w_1,...,w_n)$ is a squarefree vector in 
${\mathbb N}^n$ (i.e a vector in which $0\le w_i\le1$ for $i=1,...,n$), then we define 
its norm $|{\bf w}|$ by $|{\bf w}|=\sum_{i=1}^n w_i$. In this way, the cardinality 
$|V|$ of $V$ equals the norm of the characteristic vector $v(V)$.

Throughout the paper we denote by $R$ the polynomial ring $k[x_1,...,x_n]$ over some 
field $k$, where $n$ is the number of vertices of a hypergraph considered at the 
moment. We recall that the ring $R$ is in a natural way both $\mathbb{N}$- and 
$\mathbb{N}^n$-graded. Employing the ideas above, we may think of an edge 
$E_i$ of a hypergraph as a monomial $x^{E_i}=\prod_{j\in E_i}x_j$ in $R$. We use this 
notion to associate an ideal $I(\mc H)\sse R$ to a hypergraph $\mc H$. 
The {\bf edge ideal}, $I(\mc H)$, of a hypergraph $\mc H$ is the ideal 
$(x^{E_i} ; E_i\in\mc{E(H)})\sse R$, generated ``by the edges'' of ${\mc H}$.\\

The edge ideal was first introduced by R.~Villarreal in \cite{Vi}, in the case of simple 
graphs. Since then, edge ideals have been studied widely, see for instance 
\cite{Fa1, Fa2, VT, VT2, Vi2, Zh}. In \cite{VT} the authors give some nice recursive 
formulas for computing Betti numbers. Furthermore, their techniques illustrate 
both some obstacles that occur when you try to generalize graph theoretical results 
to hypergraph theoretical, as well as ways of getting around such obstacles. 

Another way of using hypergraphs to reveal connections between commutative algebra and 
combinatorics was introduced by S.~Faridi in \cite{Fa1}. There, Faridi consider the set of 
facets of a simplicial complex as a hypergraph. In this way a simplicial complex may be 
thought of as a ``higher dimensional'' graph. See \cite{Fa1, Fa2, Zh} for details 
and examples.\\

Recall that an {\bf (abstract) simplicial complex} on vertex set $[n]$ is a collection $\Delta$ of 
subsets of $[n]$ with the property that $F\in\Delta,\,G\sse F \Rightarrow G\in\Delta$. 
The elements of $\Delta$ are called the {\bf faces} of the complex and the maximal 
(under inclusion) faces are called {\bf facets}. The {\bf dimension} $\dim F$ of a 
face $F$ in $\Delta$ is defined to be $|F|-1$, and the dimension of $\Delta$ is defined as 
$\dim\Delta=\max\{\dim F ;\, F\in \Delta\}$. The $r$-skeleton of $\Delta$ is the 
collection of faces of dimension at most $r$. Note that the empty set $\emptyset$ 
is the unique $-1$ dimensional face of every complex that is not the void 
complex $\{\}$ which has no faces. The dimension of the void complex may be defined 
as $-\infty$.\\
The dimension $\dim_R M$ of a $R$-module $M$, is by definition the Kr\"ull dimension of 
$R/{\t{Ann}}\,M$.  

Given a simplicial complex $\Delta$, we denote by $\mc{C}.(\Delta)$ its reduced chain complex, and 
by {$\tilde{H}_n(\Delta;k)=Z_{n}(\Delta)/B_n(\Delta)$} its $n$'th reduced homology group 
with coefficients in the field $k$. In general we could use an arbitrary abelian group instead of 
$k$, but we will only consider the case when the coefficients lie in a field. For convenience, 
we define the homology of the void complex to be zero.

If $X$ and $Y$ are two sets, we denote their disjoint union by $X\sqcup Y$. Thus, suppose 
we have the two sets $[n]$ and $[m]$. They both contain the number 1, but in 
$[n]\sqcup [m]$ these two 1's are considered as distinct objects. 

Let $\Delta$ and $\Gamma$ be simplicial complexes on the disjoint vertex sets 
$\{x_1,...,x_n\}$ and $\{y_1,...,y_m\}$ respectively. We define the {\bf join} 
$\Delta\ast\Gamma$ of $\Delta$ and $\Gamma$ to be the simplicial complex on vertex set 
$\{x_1,...,x_n,y_1,...,y_m\}$ having faces 
$\{x_{i_1},...,x_{i_r},y_{j_1},...,y_{j_s}\}$, where $\{x_{i_1},...,x_{i_r}\}$ and 
$\{y_{j_1},...,y_{j_s}\}$ are faces of $\Delta $ and $\Gamma$ respectively. 

If $n\in\mathbb{N}$ we denote by $\Delta_n$ the full simplex on 
$n$ vertices. That is, the simplicial complex on $n$ vertices in which every subset 
of $[n]$ is a face. According to this we may think of the empty complex as a simplex 
on zero vertices. 

Given a simplicial complex $\Delta$ on $[n]$ and a subset $V\sse [n]$, we denote by $\Delta_V$ the 
simplicial complex on vertex set $V$, with faces 
$\{F\in\Delta ; F\sse V\}$. We call this the restriction of $\Delta$ to $V$. If ${\bf j}=
(j_1,...,j_n)$ is a squarefree vector in ${\mathbb N}^n$, by $\Delta_{\bf j}$ we mean the 
restriction to the set $V\sse [n]$ whose characteristic vector is ${\bf j}$. 

Now, let $\Delta$ be a simplicial complex on $[n]$. The 
{\bf Stanley-Reisner ring} $R/I_\Delta$ of $\Delta$ is the quotient of the ring 
$R=k[x_1,...,x_n]$ by the {\bf Stanley-Reisner ideal} 
\[
I_\Delta=(x^F ;\, F\not\in\Delta)
\]
generated by the non faces of $\Delta$. 

Let ${[n]\choose k}$ denote the set of all $k$-subsets (that is, subsets of 
cardinality $k$) of $[n]$. If $n<k$ we interpret this as being empty. Furthermore, 
we let ${n\choose k}$ denote the cardinality of ${[n]\choose k}$, so ${n\choose k}=0$ 
if $n<k$.\\

In section 2 we recall some basics that we will use throughout the paper, while 
section 3 is where the main result are found. In Theorems 3.1 and 3.5, respectively, we 
compute the Betti numbers of the $d$-complete and the $d$-complete multipartite hypergraphs, 
respectively. These results are very natural generalizations of their graph theoretical 
counterparts. By considering the independence complexes, the ideas behind the proofs 
becomes transparent. 
In section 3.4 we give a natural definition of a product on hypergraphs. This 
in turn lets us compute the Betti numbers of the $d(a_1,...,a_t)$-complete hypergraph. All 
these hypergraphs are in one way or the other a natural generalization of the ordinary 
complete graph $K_n$. In the final section, section 3.6, we define a class of hypergraphs 
that actually contain all the previously considered ones. We show that the hypergraph 
algebra, $R/I(\mc{H})$, corresponding to such hypergraph, has linear resolution.

\section{Preliminaries}

Here we recall some results and definitions which will be used
throughout the paper.

\subsection{Hypergraphs and independence complexes}
Our general reference concerning hypergraphs is Berge \cite{Be}. In this paper we
will only consider simple hypergraphs, as defined in theintroduction. Thus, hypergraph 
will always mean simple hypergraph.

Let $\mc H$ be a hypergraph. A {\bf subhypergraph} $\mc K$ of $\mc
H$ is a hypergraph such that $\mc{X}(\mc{K})\sse\mc{X(H)}$, and 
$\mc{E(K)}\sse\mc{E(H)}$. If $\mc Y\sse\mc{X}$, the {\bf induced
hypergraph on} $\mc Y$, $\mc{H}_{\mc{Y}}$, is the subhypergraph with 
$\mc{X(H_Y)}=\mc{Y}$ and with $\mc{E(H_Y)}$ consisting of the edges of $\mc{H}$ 
that lies entirely in $\mc{Y}$. A hypergraph $\mc H$ is
said to be {\bf $d$-uniform} if $|E_i|=d$ for every edge $E_i\in \mc{E(
H)}$. Note that a $2$-uniform hypergraph is just an ordinary simple
graph.\\

Let ${\mc H}=([n],\mc{E(H)})$ be a hypergraph and consider the edge
ideal $I(\mc H)\sse R$. Note that $R/I(\mc H)$ is precisely the
Stanley-Reisner ring of the simplicial complex
\[
\Delta(\mc H)=\{F\sse [n] ; E\not\sse F,\, \forall E\in \mc{E(H)}\}.
\]
This is called the {\bf independence complex} of
$\mc H$. Note that the edges in $\mc H$ are precisely the minimal
non faces in $\Delta(\mc H)$.

Let $\Delta$ be an arbitrary simplicial complex on $[n]$. We then
define the {\bf Alexander dual simplicial complex} $\Delta^\ast$ to
$\Delta$ by
\[
\Delta^\ast=\{F\sse [n] ; [n]\ssm F \not\in\Delta\}.
\]
Note that $(\Delta^\ast)^\ast=\Delta$.

\subsection{Resolutions and Betti numbers}
To every finitely generated graded module $M$ over the polynomial ring $R=k[x_1,...,x_n]$, 
we may associate a {\bf minimal ($\mathbb{N}$-)graded free resolution}
\[
0\to {\bigoplus}_j R(-j)^{\beta_{l,j}(M)}\to{\bigoplus}_j
R(-j)^{\beta_{{l-1},j}(M)} \to\cdots\to{\bigoplus}_j
R(-j)^{\beta_{0,j}(M)}\to M\to 0
\]
where $l\le n$ and $R(-j)$ is the $R$-module obtained by shifting the degrees of
$R$ by $j$. Thus, $R(-j)$ is the graded $R$-module in which 
the grade $i$ component $(R(-j))_i$ is $R_{i-j}$.\\
The natural number $\beta_{i,j}(M)$ is
called the $ij$'th $\mathbb{N}$-{\bf graded Betti number} of $M$. If $M$ is multigraded 
we may equally well consider the $\mathbb{N}^n$-graded minimal free resolution and 
Betti numbers of $M$. The difference lies just in the fact that we now use
multigraded shifts $R(-\bf j)$ instead of $\mathbb{N}$-graded ones.
The {\bf total} $i$'th Betti number is $\beta_i(M)=\sum_j
\beta_{i,j}$. For further details on resolutions, graded rings and Betti numbers, we 
refer the reader to \cite{BH}, sections 1.3 and 1.5.\\
The {\bf projective dimension} $\t{pd}(M)$ of $M$ is 
$\t{pd}(M)=\max\{i ; \exists\, \beta_{i,j}(M)\neq 0\}$.\\
\\
The Betti numbers of $M$ occur as the dimensions of certain vector
spaces over $k=R/m$, where $m$ is the unique maximal graded ideal in
$R$. Accordingly, the Betti numbers (and then of course the
projective dimension) in general depend on the characteristic of
$k$.\\
A minimal free resolution of $M$ is said to be {\bf linear} if
for $i>0$, $\beta_{i,j}(M)=0$ whenever
$j\neq i+d-1$ for some fixed natural number $d\ge1$. In this paper we only consider resolutions 
of quotient rings $R/I$. Hence, the interesting parts of the resolutions are the degrees 
greater than zero. In the variuos formulas for Betti numbers that we give, we thus assume 
that $i>0$.

In connection to this we mention the {\it Eagon-Reiner theorem}.

\begin{thm} Let $\Delta$ be a simplicial complex and $\Delta^\ast$ its Alexander dual 
complex. Then $R/I_\Delta$ is Cohen-Macaulay if and only if 
$R/I_{\Delta^\ast}$ has linear minimal free resolution.
\end{thm}
\begin{proof}
See \cite{ER}, Theorem 3.
\end{proof}

Since there is a 1-1 correspondence between Stanley-Reisner rings
(or equivalently squarefree monomial ideals) and simplicial
complexes, we get a 1-1 correspondence between simple hypergraphs
and Stanley-Reisner rings as well. This enables us to talk about resolutions, 
Betti numbers, and projective dimensions of hypergraphs.\\
By a resolution, a Betti number, or the projective dimension of a hypergraph $\mc H$, we
mean ditto of  $R/I(\mc H)$. Thus
$\beta_{i,j}(\mc H)=\beta_{i,j}(R/I(\mc H))$ and $\t{pd}(\mc
H)=\t{pd}(R/I(\mc H))$.

 One further result which we will use later on is the {\it
Auslander-Buchsbaum formula}. If $R$ is a finitely generated graded $k$-algebra for 
some field $k$ and $M\neq0$ a finitely generated graded $R$-module with $\t{pd}(M)<\infty$, 
then the formula asserts that
\[
\t{pd}(M)+\t{depth}(M)=\t{depth} R.
\]
For a proof, see \cite{BH}, Theorem 1.3.3.

\subsection{Hochster's formula}
In topology one defines Betti numbers in a somewhat different
manner. {\it Hochster's formula} provides a link between these and the
Betti numbers defined above. Hochster's formula will turn out to be a
very useful tool of ours.

\begin{thm}{(Hochster's formula).} Let $R/I_\Delta$ be the Stanley-Reisner ring of a
simplicial complex $\Delta$. The non-zero Betti numbers of
$R/I_\Delta$ are only in squarefree degrees $\bf j$ and may be
expressed as
\[
\beta_{i,{\bf j}}(R/I_\Delta)=\dim_k\tilde{H}_{|{\bf
j}|-i-1}(\Delta_{\bf j};k).
\]
Hence the total $i$'th Betti number may be expressed as
\[
\beta_i(R/I_\Delta)=\sum_{V\sse[n]}\dim\tilde{H}_{|V|-i-1}(\Delta_V;k).
\]
\end{thm}
\begin{proof}
See \cite{BH}, Theorem 5.5.1.
\end{proof}
If one has $\mathbb{N}^n$-graded Betti numbers, it is easy to obtain
the $\mathbb{N}$-graded ones via
\[
\beta_{i,j}(R/I_\Delta)=\sum_{\substack{ {\bf j'} \in {\mathbb{N}^n} \\
|{\bf j'}|=j}} \beta_{i,{\bf j'}}(R/I_\Delta).
\]
Thus,
\[
\beta_{i,j}(R/I_\Delta)=\sum_{\substack{V\sse [n]\\
|V|=j}}\dim\tilde{H}_{|V|-i-1}(\Delta_V;k).
\]

\subsection{The Mayer-Vietoris sequence}

Recall that if we have an exact sequence of complexes,\footnote{That is, complexes of 
modules over some ring $R$.}
\[
{\bf 0}\to{\bf L}\to{\bf M}\to{\bf N}\to{\bf 0}
\]
there is a long exact (reduced) homology sequence associated to it
\[
\cdots\to H_r(N)\to H_{r-1}(L)\to H_{r-1}(M)\to H_{r-1}(N)\to\cdots .
\]
Later in this paper we will have great use of this homology sequence in the special case 
where it is associated to a simplicial complex as follows.\\

Suppose we have a simplicial complex $N$ and two subcomplexes $L$ and $M$, such that $N=L\cup M$. 
This gives us an exact sequence of (reduced) chain complexes
\[
0\to\mc{C}.(L\cap M)\to\mc{C}.(L)\oplus\mc{C}.(M)\to\mc{C}.(N)\to0.
\]
The non trivial maps here are defined by $x\mapsto (x,-x)$ and $(x,y)\mapsto x+y$.

The long exact (reduced) homology sequence associated to this particular sequence, is called the 
Mayer-Vietoris sequence. The reason that we will have great use of the Mayer-Vietoris sequence 
is that in the cases that we will consider, almost always some of the considered chain complexes  
will turn out to be very easy to handle. More about 
the Mayer-Vietoris sequence can be found in \cite{Ma}, section 4.4.

\subsection{K\"unneth's tensor formula}
If complexes $\bf{L}$ and $\bf{M}$ are given, then the tensor product $\bf{L}\otimes\bf{M}$
may be constructed and given the structure of a complex as well. The degree $n$ component is 
defined as $({\bf L}\otimes{\bf M})_n=\sum_{r+s=n}L_r\otimes M_s$. Now, suppose that we 
are considering chain complexes corresponding to simplicial complexes $L$ and $M$. 
It is a natural question to ask if the (reduced) homology of the tensor product 
$\mc{C}.(L)\otimes\mc{C}.(M)$ in some way is related to the (reduced) homologies of $L$ and $M$. 
The answer is given by {\it K\"unneth's tensor formula} (\cite {Ma} Theorem 10.1), which 
under suitable\footnote{For example when the coefficients of the homology groups are in a 
field $k$.} circumstances says that
\[
\tilde{H}_n(\mc{C}.(L)\otimes\mc{C}.(M))=\bigoplus_{\substack{r+s=n\\r,s\ge0}}
\tilde{H}_r(L)\otimes\tilde{H}_s(M).
\]
We will use of this formula in connection to the join operation. It is easy to verify 
that the chaincomplex $\mc{C}.(L\ast M)$ of the join of two simplicial complexes $L$ and $M$, 
is isomorphic to the tensor product $(\mc{C}.(L)\otimes\mc{C}.(M))(-1)$. This is the same as 
the complex $(\mc{C}.(L)\otimes\mc{C}.(M))$ if we just shift the degree by 1.

\subsection{Some results on induced hypergraphs}
The formulas we have encountered so far actually yield
a couple of easy results.

Let $\mc H$ be a $d$-uniform hypergraph. We say that two edges $E$
and $E'$ are disjoint if $E\cap E'=\emptyset$. Then, by considering
the Taylor resolution (see \cite{BPS}) of
$R/I(\mc H)$, one can prove the following results, which are essentially due 
to Jacques.

\begin{prop} Let $\mc H$ be a $d$-uniform hypergraph. Then 
$\beta_{i,id}(\mc H)$ equals the number of induced hypergraphs that 
consist of $i$ disjoint edges.
\end{prop}
\begin{proof}
For $d=2$ this is Theorem 3.3.5 in \cite{Ja}. The proof given there holds also for $d>2$.
\end{proof}

\begin{prop}Let $\mc H=([n],\mc{E(H)})$ be a hypergraph and $\mc K$ an
induced hypergraph. Then
\[
\beta_{i,j}(\mc K)\le\beta_{i,j}(\mc H).
\]
\end{prop}
\begin{proof}
Since $\mc{K}=\mc{H_Y}$ for some $\mc{Y}\sse[n]$, we have
\[
\beta_{i,j}(\mc
H)=\sum_{\substack{V\sse[n]\\|V|=j}}\dim_k\tilde{H}_{|V|-i-1}(\Delta(\mc H)_V;k)
\ge\sum_{\substack{V\sse\mc{Y}\\|V|=j}}
\dim_k\tilde{H}_{|V|-i-1}(\Delta(\mc K)_V;k)=\beta_{i,j}(\mc K).
\]
\end{proof}

\begin{cor} Let $\mc H=([n],\mc{E(H)})$ be a hypergraph and $\mc K$ an
induced hypergraph. Then
\begin{eqnarray*}
\beta_i(\mc K)&\le&\beta_i(\mc H)\\
\mathrm{pd }(\mc K)&\le& \mathrm{pd }(\mc H).
\end{eqnarray*}
\end{cor}

\section{Various complete hypergraphs}
In \cite{Ja} Jacques obtains nice descriptions of the Betti numbers of some 
special families of graphs. We will generalize some of these to hypergraph analogues.

\subsection{The $d$-complete hypergraph}
The {\bf complete graph} $K_n$ on $n$ vertices is a familiar object
to all who have encountered at least some graph theory. Since an
ordinary simple graph is $2$-uniform, it seems reasonable to
consider $d$-uniform hypergraphs when
seeking a hypergraph counterpart.\\
We make the following definition. The {\bf $d$-complete} hypergraph $K_n^d$ on $n$ 
vertices is the $d$-uniform hypergraph with $\mc{E}(K_n^d)={[n]\choose d}$. 
We will now compute the Betti numbers of $K_n^d$. 

\begin{thm} The $\mathbb{N}$-graded Betti numbers of the $d$-complete
hypergraph $K_n^d$ on $n$ vertices are independent of the
characteristic of the field $k$ and may be written as
\[
\beta_{i,j}(K_n^d)=\left\{\begin{array}{lll}{n\choose
j}{{j-1}\choose
{d-1}} &if& j=i+(d-1)\\
\\
0 &if& j\neq i+(d-1).
\end{array}\right. 
\]
\end{thm}
\begin{proof}
Hochster's formula says
\[
\beta_{i,j}(K_n^d)=\sum_{\substack{V\sse [n]\\
|V|=j}}\dim_k\tilde{H}_{|V|-i-1}(\Delta(K_n^d)_V;k).
\]
It follows from the definitions that $\Delta(K_n^d)$ is the
$(d-2)$-skeleton of $\Delta_n$. In the same way, $\Delta(K_n^d)_V$ is the
$(d-2)$-skeleton of $(\Delta_n)_V\cong\Delta_{|V|}$. Thus, the complexes $\Delta(K_n^d)_V$ 
can only have non zero homology in degrees less than or equal to $d-2$. 
But, by considering a minimal resolution of $K_n^d$, it is also clear that 
$\beta_{i,j}(K_n^d)=0$ if $j<i+(d-1)$. This is simply because the 
generators of $I(K_n^d)$ have degree $d$. Hence, we have a linear resolution,
and $\beta_{i,j}(K_n^d)\neq0$ only if $j=i+(d-1)$.

Now, consider $\Delta(K_n^d)_V$ for some $V\sse [n]$.
It is clear that every cycle in $Z_{d-2}(\Delta(K_n^d))$ is a linear combination 
of ``elementary cycles'', by which we mean the derivatives of $(d-1)$-simplices in
$(\Delta_n)_V$. Denote this generating set by ${\mc G}_V$. 

We note that we may actually extract a smaller generating set out of $\mc{G}_V$. Namely, 
we claim that it is enough to consider the elements that contain a fixed 
vertex $x\in V$ (by containing $x$ we mean that some term in the cycle contains $x$). 
Denote this set by $\mc{G}_V(x)$ and consider an element $\partial(\{x_1,...,x_d\})$ 
in $\mc{G}_V$, that do not contain $x$. This cycle is a linear 
combination of elements in $\mc{G}_V(x)$, which may be seen by first forming the cone 
$x\ast\{x_1,...,x_d\}$, and then taking the derivative 
of the $(d-1)$-skeleton of this cone. This proves our claim.

Furthermore, we may easily show that the images $\bar{\sigma}$ in the homology group 
$\tilde{H}_{d-2}(\Delta(K_n^d);k)$ of the elements $\sigma\in\mc{G}_V(x)$ are linearly 
idependent. Assume that $\sum_{i=1}^t a_i\bar{\sigma_i}=0$, $a_i\in k=R/m$ (where $m$ 
is the unique graded maximal ideal of $R$) and $\sigma_i\in {\mc G}_V(x)$. Every $\sigma_i$ 
contains a unique term which does not contain $x$. This is because 
$\sigma_i=\partial(\Sigma_i)$, where $\Sigma_i$ is a $(d-1)$ simplex. Hence 
$a_i=0$ for every $i=1,...,t$. 

Now we are done since the cardinality of $\mc{G}_V(x)$ clearly is 
${{j-1}\choose{d-1}}$, and the number of $j$-sets $V$ are $n\choose j$.

\end{proof}

Due to Corollary 3.2, the above result also follows from Theorem 1 in \cite{He1}. 
Corollary 3.2 seems to be well known, but we did not manage to find a previously 
published proof.\\

Since $\Delta(K_n^d)$ has a specially nice structure, it is easy to determine its 
Alexander dual. As the minimal non-faces of $\Delta(K_n^d)$ are all 
$\{x_{i_{1}},...,x_{i_{d}}\}$, $x_{i_{j}}\in[n]$, the facets of $\Delta(K_n^d)^\ast$ 
are all $\{x_{i_{1}},...,x_{i_{n-d}}\}$, 
$x_{i_{j}}\in [n]$. Whence $\Delta(K_n^d)^\ast\cong\Delta(K_n^{n-d+1})$.

\begin{cor} The ring $R/I(K_n^d)$ is 
Cohen-Macaulay and we have
\[
\beta_i(K_n^d)={n\choose j}{{j-1}\choose {d-1}}
\]
\[
\mathrm{pd}(K_n^d)=n-(d-1)
\]
where $j=i+(d-1)$.
\end{cor}
\begin{proof}
The last two claims follows directly from the theorem. We know, by the 
Eagon-Reiner theorem, that a Stanley-Reisner ring 
$R/I_\Delta$ of a simplicial complex $\Delta$ has a linear resolution precisely 
when the Stanley-Reisner ring $R/I_{\Delta^\ast}$ of the Alexander dual complex 
is Cohen-Macaulay. Since $\Delta(K_n^d)^\ast\cong\Delta(K_n^{n-d+1})$ we are done.

\end{proof}

One should note that $\Delta(K_n^d)$ is in fact shellable. A shelling is easy to construct 
using the lexicographic order on $n$-tuples.

\begin{cor} The ring $R/I_{\Delta(K_n^d)^\ast}$ is Cohen-Macualay and we have
\begin{eqnarray*}\dim \Delta(K_n^d)^\ast&=&n-d-1\\
\dim_R(R/I_{\Delta(K_n^d)^\ast})&=&n-d\\
\mathrm{pd}(R/I_{\Delta(K_n^d)^\ast})&=&d.
\end{eqnarray*}
\end{cor}

\begin{proof}
The Cohen-Macaulayness is now clear and the first equation follows from the 
definitions. The second equation follows from the first one since $\dim_R R/I_\Delta=
\dim\Delta +1$ for any simplicial complex $\Delta$ (see \cite{BH}, Theorem 5.1.4). 
The second equation and the Cohen-Macaulayness together imply the third equation.
\end{proof}

In \cite{Ja} Jacques studies the graph algebra of $K_n$, which we denote $K_n^2$, and 
obtains the formula
\[
\beta_{i,j}(K_n)=\left\{\begin{array}{lll}{n\choose j}i &if&j=i+1\\
\\
0 &if& j\neq i+1.
\end{array}\right.
\]
Note that this is a special case of our formula for $\beta_{i,j}(K_n^d)$; just put $d=2$ 
and use the fact that $j=i+(d-1)$.

\subsection{The $d$-complete multipartite hypergraph}
Perhaps almost as familiar as the complete graph $K_n$, is the 
{\bf complete multipartite graph} $K_{n_1,...,n_t}$ on a vertex set which is a disjoint union of 
$t$ sets $[n_i]$ , with cardinality $n_i$, respectively. Contrary to the situation of the 
complete graph, it is not clear how to generalize to hypergraphs. Again, it seems 
reasonable to look for a $d$-uniform hypergraph, but this can be done in several ways. 
In this paper we will consider a few.

We define the {\bf $d$-complete multipartite} hypergraph $K_{n_1,...,n_t}^d$ on vertex 
set $[n_1]\sqcup [n_2]\sqcup\cdots\sqcup [n_t]$, to be the $d$-uniform 
hypergraph whose edge set consists of all $d$-edges except those of the form 
$\{x_{i_1},...,x_{i_d}\}$ where $x_{i_j}\in [n_i]$ for all $j=1,...,d$.

\begin{lemma} The Stanley-Reisner ring $R/I(K_{n_1,...,n_t}^d)$ of the 
$d$-complete multipartite hypergraph has linear resolution, and 
$\beta_{i,j}(K_{n_1,...,n_t}^d)\neq0$ only if $j=i+(d-1)$.
\end{lemma}
\begin{proof}
Contrary to case of the $d$-complete hypergraph, this time there may very well 
exist $(d-1)$-faces $\{x_1,...,x_d\}$ in $\Delta(K_{n_1,...,n_t}^d)$, since 
$I(K_{n_1,...,n_t}^d)$ is not generated by all possible $d$-edges.\\
As in the proof of Theorem 3.1, $\beta_{i,j}(K_{n_1,...,n_t}^d)=0$ if $j<i+(d-1)$. 
Suppose $\beta_{i,j}(K_{n_1,...,n_t}^d)\neq0$ and $j>i+(d-1)$. Via Hochster's 
formula we conclude that there must then exist a 
non zero homology group $\tilde{H}_l(\Delta(K_{n_1,...,n_t}^d)_V;k)$, for some 
$V\sse [n_1]\sqcup [n_2]\sqcup\cdots\sqcup [n_t]$ and $l\ge d-1$. 

But a cycle in such a degree $l$ has to be a sum of cycles, each of which 
lies entirely inside one of the simplices $\Delta_{n_{i}}$ on vertices $[n_i]$, 
respectively, which has no homology at all. Thus, the cycle is a boundary, contrary 
to our assumptions.
\end{proof} 

From now on it will be understod that in a multipartite situation, i.e when a hypergraph 
$\mc{H}$ has some disjoint union $[n_1]\sqcup\cdots\sqcup [n_t]$ as vertex set, 
then $\Delta_{n_s}$ denotes the simplex on the $n_s$ vertices from the $[n_s]$-component 
of $\mc{X(H)}$. We now compute the Betti numbers of $K_{n_1,...,n_t}^d$. 

\begin{thm}  The $\mathbb{N}$-graded Betti numbers of the $d$-complete 
multipartite hypergraph $K_{n_1,...,n_t}^d$ on vertex set 
$[n_1]\sqcup\cdots\sqcup [n_t]$ are independent of the
characteristic of the field $k$ and may be written as
\[
\beta_{i,j}(K_{n_1,...,n_t}^d)=\left\{\begin{array}{lll}
{{N}\choose {j}}{{j-1}\choose {d-1}}-
\sum_{\substack{ (j_1,...,j_t)\in {\mathbb{N}^t }\\
j_1+\cdots +j_t=j}}[\prod_{s=1}^t {{n_s}\choose {j_s}}]
\cdot\sum_{s=1}^t {{j_s-1}\choose {d-1}} &if& j=i+(d-1)\\
\\
0 &if& j\neq i+(d-1)
\end{array}\right. 
\]
where $N=\sum_{s=1}^t n_s$.
\end{thm}

\begin{proof}
In order to get the notations as clear as possible, we prove here only the case where $t=2$. 
It will be obvious that the same proof holds also when $t>2$. For $t=2$ the formula in the 
theorem has the following form
\[
\beta_{i,j}(K_{n,m}^d)=\left\{\begin{array}{lll} {{n+m}\choose {j}}{{j-1}\choose
{d-1}}-\sum_{j_1=0}^j {n\choose {j_{1}}}{m\choose {j-j_{1}}}[{{j_1-1}\choose {d-1}}+
{{j-j_1-1}\choose {d-1}}] &if& j=i+(d-1)\\
\\
0 &if& j\neq i+(d-1).
\end{array}\right. 
\]

Our idea is to compare the terms $\tilde{H}_{|V|-i-1}(\Delta(K_{n,m}^d)_V;k)$ 
occuring in Hochster's formula with the corresponding terms 
$\tilde{H}_{|V|-i-1}(\Delta(K_{n+m}^d)_V;k)$ which we encountered when we computed 
$\beta_{i,j}(K_{n+m}^d)$.

We realize, simply because we have descriptions of the structures of the considered 
complexes, that
\[
\dim_k\tilde{H}_{|V|-i-1}(\Delta(K_{n,m}^d)_V;k)\le
\dim_k\tilde{H}_{|V|-i-1}(\Delta(K_{n+m}^d)_V;k)
\]
for every set $V\sse [n]\sqcup [m]$. The possible difference lies in the fact that 
there might very well be faces $F\in\Delta(K_{n,m}^d)$ such that $|F|\ge d$. 
This would result in a non zero boundary group $B_{d-2}(\Delta(K_{n,m}^d))$ 
in the chain complex of $\Delta(K_{n,m}^d)$. 

It is an elementary fact that
\[
\dim_k\tilde{H}_{d-2}(\Delta(K_{n,m}^d)_V;k)=\dim_k Z_{d-2}(\Delta(K_{n,m}^d)_V)-
\dim_k B_{d-2}(\Delta(K_{n,m}^d)_V).
\]
Since the cycle groups $Z_{d-2}(\Delta(K_{n+m}^d)_V)$ and 
$Z_{d-2}(\Delta(K_{n,m}^d)_V)$ 
clearly coincide and since $B_{d-2}(\Delta(K^d_{n+m})_V)=0$, we only have to compute the 
dimension over $k$ of 
$B_{d-2}(\Delta(K_{n,m}^d)_V)$.\\

If we write $V=V_1\sqcup V_2$, where $V_1\sse [n]$ and $V_2\sse [m]$, it is clear that 
\[
B_{d-2}(\Delta(K_{n,m}^d)_V)=B_{d-2}(\Delta(K_{n,m}^d)_{V_1})\oplus 
B_{d-2}(\Delta(K_{n,m}^d)_{V_2}).
\]
This is because the potential $(d-1)$-faces of 
$\Delta(K_{n,m}^d)$ lies either in $\Delta(K_{n,m}^d)_{[n]}$ or in 
$\Delta(K_{n,m}^d)_{[m]}$, which are disjoint.

Now, we have already proved how to compute $\dim_k B_{d-2}(\Delta(K_{n,m}^d)_{V_\nu})$,  
$\nu=1,2$. This was done when we computed the Betti numbers of $K_n^d$. Thus, 
\[
\dim_k B_{d-2}(\Delta(K_{n,m}^d)_{V_1})={{|V_1|-1}\choose {d-1}}
\]
\[
\dim_k B_{d-2}(\Delta(K_{n,m}^d)_{V_2})={{|V_2|-1}\choose {d-1}}.
\]
If we put $|V_1|=j_1$ the theorem follows as we simply sum over all possible 
$V\sse [n]\sqcup [m]$. 

\end{proof}

\begin{cor} Given $K_{n_1,...,n_t}^d$ with $N=\sum_{s=1}^t n_s\ge d$, we have
\[
\beta_i(K_{n_1,...,n_t}^d)={{N}\choose {j}}{{j-1}\choose {d-1}}-
\sum_{\substack{ (j_1,...,j_t)\in {\mathbb{N}^t }\\
j_1+\cdots +j_t=j}}[\prod_{s=1}^t {{n_s}\choose {j_s}}]
\cdot\sum_{s=1}^t {{j_s-1}\choose {d-1}}
\]
\[
\mathrm{pd}(K_{n_1,...,n_t}^d)=N-(d-1)
\] 
where $j=i+(d-1)$.
\end{cor}
\begin{proof}
The fact that $\t{pd}(K_{n_1,...,n_t}^d)\le N-(d-1)$ follows from directly 
from the formula. By putting $j=N$ we get
\[
{{N-1}\choose {d-1}}-\sum_{s=1}^t {{n_s-1}\choose {d-1}}.
\]
This expression is strictly greater than 0, which we may prove as follows. 
Consider the set $[n_1]\sqcup [n_2]\sqcup\cdots\sqcup [n_t]$ of $N$ elements. Pick an 
arbitrary element and remove it from the set. The first term above count the number of 
ways of choosing $d-1$ elements from the later set.\\
The sum in the above display counts the following: Start with the same set at before, 
and remove an arbitrary element $x_s$ from each one of the sets $[n_s]$. Then choose 
$(d-1)$ elements from some $[n_s]\ssm x_s$. $\sum_{s=1}^t{{n_s-1}\choose{d-1}}$ is the 
total number of such $(d-1)$-sets.\\
Clearly, the difference of these two numbers is strictly greater than 0, just consider a 
set of $(d-1)$-elements that do not lie entirely inside one set $[n_s]$. 
As we have assumed that $N\ge d$ the claim follows.   

\end{proof}
{\bf Example:}\, Denote the vertex set of $K_{2,3}^3$ by 
$\{a,b\}\sqcup \{A,B,C\}$. Then we have 
\[
\mc{E}(K_{2,3}^3)=\{abA,abB,abC,aAB,bAB,aAC,bAC,aBC,bBC\}.
\]
The Betti numbers are $\beta_0(K_{2,3}^3)=1,\,\beta_1(K_{2,3}^3)=9,\, \beta_2(K_{2,3}^3)=13,\, 
\beta_3(K_{2,3}^3)=5$.\\ 

By construction, the edges in a hypergraph $\mc{H}$ are the minimal non faces in 
$\Delta(\mc{H})$. This makes it easy to determine the facets in $\Delta(\mc{H})^\ast$. 
As one easily realizes, they are the complements of the edges. Considering this, we get the 
following expression for the Alexander dual complex.
\[
\Delta(K_{n_1,..,n_t}^d)^\ast=
\bigcup_{s=1}^t(\bigcup_{l=1}^{l_s}
[\Gamma_{n_s-(l+1)}(n_s)\ast\Gamma_{n_1+\cdots +\widehat{n_s}+\cdots +n_t-d+l-1}
(n_1,...,\widehat{n_s},...,n_t)]
\]
where $\Gamma_r(n_s)$ is the $r$-skeleton of $\Delta_{n_s}$, 
$\Gamma_r(n_1,...,\widehat{n_s},...,n_t)$ is the $r$-skeleton of 
$\Delta_{n_1}\ast\cdots\ast\widehat{\Delta_{n_s}}\ast\cdots\ast\Delta_{n_t}$,\quad
$\widehat{\cdot}$\,\, means omit and $\l_s=\min\{d-1,n_s\}$.

\begin{cor} The ring $R/I_{\Delta(K_{n_1,...,n_t}^d)^\ast}$ is Cohen-Macaulay and we have
\begin{eqnarray*} \dim\Delta(K_{n_1,...,n_t}^d)^\ast&=&N-d-1\\
\dim_R(R/I_{\Delta(K_{n_1,...,n_t}^d)^\ast})&=&N-d\\
\mathrm{pd}(R/I_{\Delta(K_{n_1,...,n_t}^d)^\ast})&=&d.
\end{eqnarray*}
\end{cor}

\begin{proof}
The Cohen-Macaulayness follows from Lemma 3.4 and the Eagon-Reiner 
theorem. By considering the above description of the Alexander dual, the first equation is 
clear and implies the second. The third equation follows since
$R/I_{\Delta(K_{n_1,...,n_t}^d)^\ast}$ is Cohen-Macaulay.

\end{proof}

Also in this case we have generalized a formula given by Jacques in \cite{Ja}. 
By studying the graph algebra of $K_{n,m}$ he obtains the formula 
\[
\beta_{i,j}(K_{n,m})=\left\{\begin{array}{lll}
\sum_{j_1=1}^{j-1}{n\choose j_1}{m\choose {j-j_1}} &if& j=i+1\\
\\
0 &if& j\neq i+1.
\end{array}\right.
\]
A priori this looks quite different from our result. But, if one put $d=2$ and 
remember that ${n\choose d}$ is defined as 0 if $n<d$, our formula simplifies 
immediately to this one.

Contrary to when we considered $\Delta(K_n^d)$, the structure of the Alexander dual 
$\Delta(K_{n_1,...,n_t}^d)^\ast$ is not transparent. One immediate question that appear 
is: When, if at all, does the Stanley-Reisner ring 
of $\Delta(K_{n_1,...,n_t}^d)$ both have linear resolution and the Cohen-Macaulay 
property? Since we already know that all considered resolutions are linear, we only 
have to think about the Cohen-Macaulay property.

\begin{lemma} Let $N=\sum_{s=1}^t n_s$ and 
$n_s\le d-1$ for $s=1,...,t$. Then $K_N^d=K_{n_1,...,n_t}^d$.
\end{lemma}
\begin{proof}
${\mc E(K_{n_1,...,n_t}^d)}={\mc E(K_N^d)}$ and 
$\mc{X}(K_{n_1,...,n_t}^d)=\mc{X}(K_N^d)$.
\end{proof}

\begin{prop} The Stanley-Reisner ring $R/I(K_{n_1,...,n_t}^d)$ of a 
$d$-complete multipartite hypergraph on vertex set 
$[n_1]\sqcup [n_2]\sqcup\cdots\sqcup [n_t]$ is Cohen-Macaulay precisely when 
$n_s\le d-1$ for all $s=1,...,t$.
\end{prop}
\begin{proof} The Auslander-Buchsbaum formula tells us that
\[
\t{pd}(K_{n_1,...,n_t}^d)+\t{depth}_R(K_{n_1,...,n_t}^d)=N
\]
where $N=\sum_{s=1}^t n_s$. Since we already have computed the projective 
dimension, the above formula says
\[
\t{depth}_R(K_{n_1,...,n_t}^d)=d-1
\]
and it is clear that $\dim\Delta(K_{n_1,...,n_t}^d)=\max\{n_i-1, d-2 ;i=1,...,t\}$. 
Thus, since $\t{depth}_R M\le\dim_R M$ holds for every finitely generated 
$R$-module $M$, $R/I(K_{n_1,...,n_t}^d)$ is Cohen-Macaulay precisely when $n_s\le d-1$ 
for all $s=1,...,t$. Furthermore, according to the lemma, we have $K_{n_1,...,n_t}^d=K_{N}^d$. 
\end{proof}

\subsection{Hilbert series}
Let $M$ be a ${\mathbb{N}}$-graded module (${\mathbb{N}^n}$-graded would work 
equally well). The {\bf Hilbert series} $H_M(t)$ measure the 
dimensions over $k=R/m$ of the graded pieces $M_i$ of $M$. More algebraically: 
Let $M$ be such that every graded piece $M_i$ has finite dimension over $k$. 
Then $H_M(t)$ is the formal power series
\[
H_M(t)=\sum_{i\in\mathbb{N}}\dim_k(M_i)t^i.
\]
The following is a well known result. See for example \cite{BH}, Theorem 4.1.13.

\begin{lemma} Let $R$ be the polynomial ring $k[x_1,...,x_n]$ over a field $k$ and 
consider a finitely generated $\mathbb{N}$-graded $R$-module $M$. 
Then
\[
H_M(t)=\frac{S_M(t)}{(1-t)^n}.
\]
where $S_M(t)=\sum_{i,j}(-1)^i\beta_{i,j}(M)t^j$.
\end{lemma}

If $M$ is the Stanley-Reisner ring of a simplicial complex $\Delta$, one may 
rather easily compute its Hilbert series. This is Corollary 1.15 in \cite{MS}.  
One gets
\[
H_{R/{I_\Delta}}(t)=\frac{1}{(1-t)^n}\sum_{r=0}^e f_{r-1}t^r (1-t)^{n-r}
\]
where $f_r$ equals the number of $r$-faces of $\Delta$ and $e=\dim\Delta+1$. 

Note that this gives a nice connection between the ``geometric'' numbers $f_r(\Delta)$ and 
the ``algebraic'' numbers $\beta_{i,j}(R/I_\Delta)$. In general though, it might be 
quite messy to handle the alternating sum of Betti numbers. But, if we consider a 
module $M$ with linear resolution, the correspondence becomes much nicer.

\begin{lemma} Let $\Delta$ be a simplicial complex such that 
$R/I_\Delta$ has a linear resolution.  Then we have
\[
\beta_{i,j}(R/I_\Delta)=\sum_{r=0}^e (-1)^{j-i-r}f_{r-1}{{n-r}\choose{j-r}}.
\]
\end{lemma}
\begin{proof}
From Lemma 3.10 we get one expression for $(1-t)^n H_{R/I_{\Delta}}(t)$, and from the 
discussion right after that lemma we get another. Just identify the coefficient of $t^j$ 
from the two expressions.  
\end{proof}
This lemma gives us an alternative way of computing the Betti numbers of $K_n^d$ and 
$K_{n_1,...,n_t}^d$. All we need is the 
${\bf f}$-{\bf vector}=$(f_{-1}, f_0, f_1,..., f_{e-1})$. In the cases considered, the 
$f$-vectors have nice and simple descriptions.

Let us begin by considering $\Delta(K_n^d)$. Since this is the $(d-2)$-skeleton 
of $\Delta_n$, we see that $\dim\Delta(K_n^d)=d-2$. Thus $e$ equals $d-1$ in this case. 
The number of $(r-1)$-faces clearly is ${n\choose r}$, so the $f$-vector is given by
\[
(1,n,{n\choose 2},...,{n\choose {d-1}}).
\]
According to the above, recalling that $j=i+(d-1)$ we get the formula
\[
\beta_{i,j}(K_n^d)=\sum_{r=0}^{d-1}(-1)^{(d-1)-r}{n\choose r}{{n-r}\choose {j-r}}.
\]
This is without a doubt correct, but looks completely different from our earlier expression. 
We obviously have
\[
\sum_{r=0}^{d-1}(-1)^{(d-1)-r}{n\choose r}{{n-r}\choose {j-r}}=
{n\choose j}{{j-1}\choose {d-1}}.
\]
This identity may also be proved in a combinatorial way, using the Principle of 
Inclusion-Exclusion. We give the main ideas here. The trick is to identify 
something that is counted by both sides of the identity. This something is described below.\\
\\
 1) Consider a set of $n$ elements. First choose $j$ elements of these, and then 
choose one of the $j$ and colour it. Then colour $d-1$ further elements chosen from the 
remaining $j-1$ elements. This can be done 
in $j{n\choose j}{{j-1}\choose {d-1}}$ ways.\\
\\
We now claim that the following process counts the same thing.\\
\\
2) Choose $d-1$ elements of the $n$-set and colour them. Choose $j-d+1$ elements out of 
the remaining $n-d+1$ elements not previously choosen. Then choose one of the $j$ elements 
choosen so far and colour it. This can be done in $j{n\choose {d-1}}{{n-d+1}\choose 
{j-d+1}}$ ways. We realize that we have counted more coloured sets than in 1) in this 
process, for example those in which only $d-1$ element became coloured. In an attempt to adjust 
this we subtract $j{n\choose {d-2}}{{n-d+2}\choose {j-d+2}}$ from 
$j{n\choose {d-1}}{{n-d+1}\choose {j-d+1}}$. This number is created using 
the same choice argument as before. Then we subtract the number of coloured sets in which 
only $d-1$ elements were coloured. But we subtract 
too much, since we also subtract the number of sets in which only $d-2$ elements is 
coloured. Thus, we have to add back. 

Continuing this process, according to the Principle of Inclusion-Exclusion, 
after a finite number of steps we will stop and the resulting number counts precisely 
the same thing as 1). Finally, we just divide every term by $j$ to obtain our identity.\\ 
The number described in 1) and 2), counts the number of ways of: Choosing a $j$-set of $[n]$ to 
form a football team, say, and then determining in how many ways one can have $d$ of the 
players on the field, one of which is to be choosen as goalkeeper.\\

Note that the above arguments makes sense only if $j\ge d$. However, according to our 
earlier investigations, this is quite natural.  

We also obtain a different formula for the Betti numbers 
$\beta_{i,j}(K_{n_1,...,n_t}^d)$. Just as before, we only need to compute the 
$f$-vector. This is sufficient since we know that $K_{n_1,...,n_t}^d$ has linear 
resolution.

The structure of $\Delta(K_{n_1,...,n_t}^d)$ is easy to understand, and it follows that
\[
f_{r-1}(\Delta(K_{n_1,...,n_t}^d))=\left\{\begin{array}{lll}
{N\choose r} &if& r\le d-1\\
\\
\sum_{s=1}^t{n_s\choose r} &if& r\ge d
\end{array}\right.
\] 
where $N=\sum_{s=1}^t n_s$. Using the lemma and remembering that $j=i+(d-1)$, we obtain the 
following formula 
\[
\beta_{i,j}(K_{n_1,...,n_t}^d)=\sum_{r=0}^{d-1}(-1)^{(d-1)-r}{{N}\choose r}
{{N-r}\choose {j-r}}-\sum_{r=d}^{e}(-1)^{d-r}{{N-r}\choose {j-r}}
[\sum_{s=1}^t{n_s\choose r}]
\]
where $e=\max\{n_s-1, d-2 ; s=1,...,t\}$ is the dimension of $\Delta(K_{n_1,...,n_t}^d)$. 
We immediately note one thing. The first sum in the display actually is nothing but 
$\beta_{i,j}(K_N^d)={N\choose j}{{j-1}\choose {d-1}}$, where $N=\sum_{s=1}^t n_s$. 
Thus, we realize that that the second sum gives us an alternative expression for the difference 
$\beta_{i,j}(K_N^d)-\beta_{i,j}(K_{n_1,...,n_t}^d)$. 
In other words, we have an identity
\[
\sum_{r=d}^{e}(-1)^{d-r}{{N-r}\choose {j-r}}
[\sum_{s=1}^t{n_s\choose r}]=
\sum_{\substack{ (j_1,...,j_t)\in {\mathbb{N}^t }\\
j_1+\cdots +j_t=j}}[\prod_{s=1}^t {{n_s}\choose {j_s}}]
\cdot\sum_{s=1}^t {{j_s-1}\choose {d-1}}.
\]

\subsection{The Alexander dual of a join}
It is known, and proved in for example \cite{F}, that the join $\Delta\ast\Gamma$ of two 
simplicial complexes $\Delta$ and $\Gamma$ is Cohen-Macaulay precisely when both 
$\Delta$ and $\Gamma$ are Cohen-Macaulay. In that case, remember that the Eagon-Reiner 
theorem tells us that the Alexander dual compex $(\Delta\ast\Gamma)^\ast$ has linear resolution.

In this paper we consider several classes of hypergraphs with linear resolutions. Therefore, 
it would be nice to be able to describe the Alexander dual of a join since we then rather 
easily can construct more hypergraphs with linear resolutions. In this section we derive a 
description of the Alexander dual of a join, and also give a formula for the Betti numbers.\\

Let $\Delta$ and $\Gamma$ be simplicial complexes on $[n]$ and $[m]$, respectively. 
Denote the minimal non faces of $\Delta$ ($\Gamma$) by $f_i$, $i=1,...,r$ 
($g_j$, $j=1,...,s$, respectively). Remember that according to the identifications that are made in 
the introduction, we may consider the $f_i$'s ($g_j$'s) as squarefree monomials in 
$k[x_1,...,x_n]$ ($k[y_1,...,y_m]$). Using this identification, we consider the 
Stanley-Reisner ideal $I_{\Delta}\sse k[x_1,...,x_n]$ ($I_{\Gamma}\sse k[y_1,...,y_m]$) 
of $\Delta$ ($\Gamma$, respectively). 
It is well known (\cite{MS}, Theorem 1.7) that 
\[
I_{\Delta}=(f_i\,;\,i=1,...,r)=\bigcap_{f\in\Delta}{\mathfrak{m}}^{\bar{f}}
\]
where for a subset $V\sse[n]$, ${\mathfrak{m}}^V$ is the ideal 
$(x_i\,;\,i\in V)\sse k[x_1,...,x_n]$, and by $\bar{f}$ we mean $[n]\ssm f$. It is easy to 
realize that it is enough to take the intersection where $f$ is a facet of $\Delta$. 
If we consider the Alexander dual $\Delta^\ast$ in the same way, we get
\[
I_{\Delta^\ast}=(f'_i\,;\,i=1,...,r')=\bigcap_{f'\in\Delta^\ast}\mathfrak{m}^{\bar{f'}}
\]
where $f'_i$, $i=1,...,r'$ is the set of minimal non faces of $\Delta^\ast$ (analogously 
we denote by $g'_j$, $j=1,...,s'$, the minimal non faces of $\Gamma^\ast$). Note that this shows
the algebraic version of Alexander duality. The association $\Delta\mapsto\Delta^\ast$ is 
by the above equivalent to
\[
I_{\Delta}=(f_i\,;\,i=1,...,r)\mapsto\bigcap_{i=1}^r\mathfrak{m}^{f_i}=I_{\Delta^\ast}.
\]

If we consider $I_{\Delta}$ and $I_{\Gamma}$ as ideals in $k[x_1,...,x_n,y_1,...,y_m]$, it 
follows that
\[
I_{\Delta\ast\Gamma}=I_{\Delta}+I_{\Gamma}=(f_i,g_j\,;\,i=1,...,r,\,j=1,...,s).
\]
Hence, the Stanley-Reisner ideal of $(\Delta\ast\Gamma)^\ast$ is 
\[
\bigl(\bigcap_{i=1}^r\mathfrak{m}^{f_i}\bigr)\bigcap\bigl(\bigcap_{j=1}^s
\mathfrak{m}^{g_j}\bigr)=I_{\Delta^\ast}\cap I_{\Gamma^\ast}=I_{\Delta^\ast}I_{\Gamma^\ast}.
\]

So, we conclude that
\[
I_{(\Delta\ast\Gamma)^\ast}=(f'_ig'_j\,;\,i=1,...,r',\,j=1,...,s').
\]
By considering the minimal nonfaces $f'_ig'_j$ of $(\Delta\ast\Gamma)^\ast$, we realize that
\[
(\Delta\ast\Gamma)^\ast=\bigl(\Delta^\ast\ast\Delta_m\bigr)\bigcup
\bigl(\Delta_n\ast\Gamma^\ast\bigr).
\]
Note that we could have reached these conclusions also by considering the 
minimal non faces of $\Delta\ast\Gamma$. The form of the Stanley-Reisner ideal of 
$(\Delta\ast\Gamma)^\ast$ is particularly nice since the generators correspond to 
edges in certain hypergraphs.\\

Suppose that hypergraphs $\mc{H}=([n],\mc{E(H)})$ and $\mc{K}=([m],\mc{E(K)})$ are given. We 
define the {\bf product} $\mc{H\cdot K}$ of $\mc{H}$ and $\mc{K}$ to be the hypergraph 
on vertex set $[n]\sqcup [m]$ and with edges $\{x_1,...,x_r,y_1,...,y_s\}$, where 
$\{x_1,...,x_r\}$ is an edge in $\mc{H}$ and $\{y_1,...,y_s\}$ is an edge in $\mc{K}$. 
In other words, $\mc{E(H\cdot K)}$ may be thought of as the cartesian product 
$\mc{E(H)}\times\mc{E(K)}$.\\

Using the above results, we may easily prove the following theorem.

\begin{thm} Let $\mc{H}=([n],\mc{E(H)})$ and $\mc{K}=([m],\mc{E(K)})$ be $d$- and $d'$-uniform 
hypergraphs respectively. Then $\mc{H}\cdot\mc{K}$ is a $(d+d')$-uniform hypergraph, and has 
linear resolution if and only if both $\mc{H}$ and $\mc{K}$ have linear resolutions.       
\end{thm}
\begin{proof}
The fact that $\mc{H\cdot K}$ is $(d+d')$-uniform is clear from the definition. If we put 
$\Delta=\Delta(\mc{H})^\ast$ and $\Gamma=\Delta(\mc{K})^\ast$ in the results deduced just before 
the theorem, we get that
\[
\Delta(\mc{H\cdot K})=(\Delta(\mc{H})^\ast\ast\Delta(\mc{K})^\ast)^\ast.
\]
This is clear considering the minimal non faces of both sides of the equation. By the 
Eagon-Reiner theorem $(\Delta(\mc{H})^\ast\ast\Delta(\mc{K})^\ast)^\ast$ has 
linear resolution precisely when both $\Delta(\mc{H})^\ast$ and $\Delta(\mc{K})^\ast$ 
are Cohen-Macaulay. This is, again by the Eagon-Reiner theorem, the same thing as saying that 
both $\Delta(\mc{H})$ and $\Delta(\mc{K})$ have linear resolutions.

\end{proof} 

Note that the topological information in the above theorem says that
\[
\Delta(\mc{H\cdot K})=\bigl(\Delta(\mc{H})\ast\Delta_m\bigr)\bigcup
\bigl(\Delta_n\ast\Delta(\mc{K})\bigr).
\]
Now, let $V=V_1\sqcup V_2\sse [n]\sqcup [m]$ and consider the exact sequence
\[
0\to\mc{C}.((\Delta(\mc{H})\ast\Delta(\mc{K}))_V)\to
\mc{C}.((\Delta(\mc{H})\ast\Delta_m)_V)\oplus\mc{C}.((\Delta_n\ast\Delta(\mc{K}))_V)\to
\mc{C}.(\Delta(\mc{H\cdot K})_V)\to0.
\]
If $V_1$ or $V_2$ is empty, then $\Delta(\mc{H\cdot K})_V$ will not have any non zero homology. 
This is simply because there are no non faces (consider the relations in the Stanley-Reisner ring). 
Our aim is to compute the Betti numbers via Hochster's formula and hence, it is enough to consider 
the sets $V=V_1\sqcup V_2\sse [n]\sqcup [m]$ for which $V_1\cap [n]$ and $V_2\cap [m]$ both are 
non empty. But, in this case both $(\Delta(\mc{H})\ast\Delta_m)_V$ and 
$(\Delta_n\ast\Delta(\mc{K}))_V$ are cones and accordingly have no homology at all. Thus, 
if we consider the Mayer-Vietoris sequence obtained from the above exact sequence, we get that 
the following equation holds for every $V\sse [n]\sqcup [m]$, $V\cap [n]\neq\emptyset$, 
$V\cap [m]\neq\emptyset$:
\[
\tilde{H}_{r}(\Delta(\mc{H\cdot K})_V;k)\cong
\tilde{H}_{r-1}((\Delta(\mc{H})\ast\Delta(\mc{K}))_V;k).
\]
Using the results in section 2.5, it follows that
\[
\tilde{H}_{r}(\Delta(\mc{H\cdot K})_V;k)\cong
\bigoplus_{\substack{r_1+r_2=r-2\\r_1,r_2\ge0}}
\tilde{H}_{r_1}(\Delta(\mc{H})_V;k)\otimes\tilde{H}_{r_2}(\Delta(\mc{K})_V;k).
\]
Thus, by Hochster's formula, we get
\[
\beta_{i,j}(\mc{H\cdot K})=\sum_{\substack{|V|=j\\V=V_1\sqcup V_2}}
\sum_{r_1+r_2=j-i-3}\dim_k\tilde{H}_{r_1}(\Delta(\mc{H})_{V_1};k)\cdot
\dim_k\tilde{H}_{r_2}(\Delta(\mc{K})_{V_2};k).
\]
Of course, we want to extend this to products of more than two hypergraphs. This we do inductively. 

\begin{thm} The $ij$'th $\mathbb{N}$-graded Betti number of the product $\mc{H}_1\cdots\mc{H}_t$ 
of hypergraphs $\mc{H}_i$, $i=1,...,t$ on vertex sets $[n_i]$ respectively, is given by the 
following expression.
\[
\beta_{i,j}(\mc{H}_1\cdots\mc{H}_t)=\sum_{\substack{|V|=j\\V=V_1\sqcup\cdots\sqcup V_t}}
\sum_{\substack{r_1+\cdots +r_t=j-i-(2t-1)\\r_i\ge0}}
\prod_{l=1}^t\dim_k\tilde{H}_{r_l}(\Delta(\mc{H}_l)_{V_l};k).
\]
\end{thm}
\begin{proof}
We have already seen that the formula holds for $t=2$. It follows easily by induction that 
\[
\dim_k\tilde{H}_s(\Delta(\mc{H}_1\cdots\mc{H}_t)_V;k)=
\sum_{r_1+\cdots+r_t=s-2(t-1)}\prod_{l=1}^t\dim_k\tilde{H}_{r_l}
(\Delta(\mc{H}_l)_{V_l};k).
\]
Now consider the following, which by the case $t=2$ clearly holds.
\[
\beta_{i,j}(\mc{H}_1\cdots\mc{H}_{t+1})=\sum_{\substack{|V|=j\\V=V_1\sqcup\cdots\sqcup V_t}}
\sum_{s+r_{t+1}=j-i-3}\dim_k\tilde{H}_s(\Delta(\mc{H}_1\cdots\mc{H}_t)_V;k)\cdot
\dim_k\tilde{H}_{r_{t+1}}(\Delta(\mc{H}_{t+1})_V;k).
\]
By putting the expression for $\dim_k\tilde{H}_{s}(\Delta(\mc{H}_1\cdots\mc{H}_{t})_V;k)$ 
in the above formula, we easily see that the two equations 
\[
r_1+\cdots+r_t=s-2(t-1)
\]
\[
s+r_{t+1}=j-i-3
\]
may be collected into the single equation
\[
r_1+\cdots+r_{t+1}=j-i-(2(t+1)-1).
\]
By induction we are done.

\end{proof}

The above formula for the Betti numbers becomes much nicer if we know that each $\mc{H}_i$ has 
linear resolution. The effect of this is that the inner summation symbol becomes superfluous, 
this since we already know that in this case 
$\dim_k\tilde{H}_{r_l}(\Delta(\mc{H}_l)_{V_l};k)$ can only be non zero in one 
specific degree for each $l$. These pieces of information yield the degree in which 
$\dim_k\tilde{H}_{s}(\Delta(\mc{H}_1\cdots\mc{H}_{t})_V;k)$ is non zero. But, this degree 
is also expressed by the equation $r_1+\cdots+r_t=j-i-(2t-1)$. Thus, we may indeed remove the 
summation symbol and we have 

\begin{thm} Let hypergraphs $\mc{H}_i$, $i=1,...,t$, on vertex sets $[n_i]$ respectively be given.
Assume that for $i=1,...,t$, the hypergraph $\mc{H}_i$ is $a_i$-uniform with linear resolution. 
Then the $ij$'th $\mathbb{N}$-graded Betti number of the product 
$\mc{H}_1\cdots\mc{H}_t$ is given by the following expression.
\[
\beta_{i,j}(\mc{H}_1\cdots\mc{H}_t)=\sum_{\substack{|V|=j\\V=V_1\sqcup\cdots\sqcup V_t}}
\prod_{l=1}^t\dim_k\tilde{H}_{a_l-2}(\Delta(\mc{H}_l)_{V_l};k).
\]
\end{thm}
\begin{flushright}
$\Box$
\end{flushright}

\subsection{The $d(a_1,...,a_t)$-complete multipartite hypergraph}
As we mentioned before, there are many ways of generalizing the multipartite graph 
$K_{n_1,...,n_t}$ to a hypergraph analogue. We have already discussed the $d$-complete 
multipartite hypergraph $K_{n_1,...,n_t}^d$, and will now move on to consider another 
class of hypergraphs.\\

The edge set $\mc E(K_{n_1,...,n_t}^d)$ of $K_{n_1,...,n_t}^d$ consists of all $d$-edges 
except those of the form $\{x_{i_1},...,x_{i_d}\}$, $x_{i_j}\in [n_s]$ for some 
$s=1,...,t$. In the case of the ordinary graph $K_{n,m}$, this just tells us that we have 
all edges between the disjoint sets $[n]$ and $[m]$ of vertices. This one may think of as 
an edge being a choice of two vertices, prescribing a certain number of vertices in each one 
of the sets $[n]$ and $[m]$, namely one in each. This is the idea behind what we now define. 
The {\bf $d(a_1,...,a_t)$-complete multipartite} hypergraph 
$K_{n_1,...,n_t}^{d(a_1,...,a_t)}$ is the $d$-uniform hypergraph on vertex set 
$[n_1]\sqcup [n_2]\sqcup\cdots\sqcup [n_t]$ and edge set 
$\mc{E}(K_{n_1,...,n_t}^{d(a_1,...,a_t)})$ consisting of all $d$-edges such that 
precisely $a_s$ elements comes from $[n_s]$, $a_s\in\mathbb{N}$,  $a_s\ge1$, 
$\sum_{s=1}^t a_s =d$.\footnote{The $d$ occuring in the superscript in the symbol 
$K_{n_1,...,n_t}^{d(a_1,...,a_t)}$ does not have any real purpose here. However, when we 
continue our work, the $d$ will be useful since the notations will become more unified.} 
 
Proceeding in the same spirit as before we begin our investigation by showing that 
$R/I(K_{n_1,...,n_t}^{d(a_1,...,a_t)})$ has a linear resolution. First let us simplify the 
notation a bit. In what follows, $a_1,...,a_t={\bf a}$,
$n_1,...,n_t = {\bf n}$ and $d=\sum_{s=1}^t a_s$. Thus, $d(a_1,...,a_t)=d({\bf a})$ and 
$K_{\bf {n}}^{d({\bf a})}=K_{n_1,...,n_t}^{d(a_1,...,a_t)}$.

\begin{lemma} The Stanley-Reisner ring $R/I(K_{\bf n}^{d(\bf a)})$ of the 
$d({\bf a})$-complete mutlipartite hypergraph has linear resolution, and 
$\beta_{i,j}(K_{\bf n}^{d(\bf a)})\neq0$ only if $j=i+(d-1)$.
\end{lemma}
\begin{proof}
This will be proved in greater generality in section 3.6. However we give a short 
proof here as well. This is since it contains some interesting information which we will 
not get out of the more general proof in section 3.6.\\

By considering the definition of the Alexander dual complex, we immediately get the following 
expression.
\[
\Delta(K_{\bf n}^{d(\bf a)})^\ast=
\Gamma_{n_1-a_1-1}(n_1)\ast\cdots\ast\Gamma_{n_t-a_t-1}(n_t).
\]
Thus, $\Delta(K_{\bf n}^{d(\bf a)})^\ast$ is Cohen-Macaulay since we know that each 
$\Gamma_{n_s-a_s-1}(n_s)$ is Cohen-Macaulay. Now, our lemma follows by the Eagon-Reiner theorem.
\end{proof}  

We now compute the Betti numbers of $K_{\bf n}^{d(\bf a)}$. As one easily realizes, 
either from the above lemma or directly from the definition, we have that 
\[
K_{\bf n}^{d({\bf a})}=\prod_{s=1}^t K_{n_s}^{a_s}.
\]
Thus, we may apply the results from the previous section.

\begin{thm} The $\mathbb{N}$-graded Betti numbers of 
the $d(\bf a)$-complete multipartite hypergraph $K_{\bf n}^{d(\bf a)}$ on vertex set 
$[n_1]\sqcup\cdots\sqcup [n_t]$ are independent of the characteristic of the field $k$ and 
may be written as
\[
\beta_{i,j}(K_{\bf n}^{d(\bf a)})=
\left\{\begin{array}{lll}\sum_{\substack{r_1+\cdots +r_t=i+t-1\\
r_i\ge1}}[\prod_{l=1}^t \beta_{r_l, r_l+a_l-1}(K_{n_l}^{a_l})] &if& j=i+(d-1)\\
0 &if& j\neq i+(d-1).
\end{array}\right.
\]
\end{thm}
\begin{proof}
We know that $\dim_k\tilde{H}_{r_l}(\Delta(K_{n_l}^{a_l})_{V_l};k)\neq0$ only when 
$r_l=a_l-2$, and in this case, we have
\[
\dim_k\tilde{H}_{r_l}(\Delta(K_{n_l}^{a_l})_{V_l};k)={{j_l-1}\choose{a_l-1}}
\]
where $j_l=|V_l|$. Using this, the expression
\[
\sum_{\substack{|V|=j\\V=V_1\sqcup\cdots\sqcup V_t}}
\sum_{\substack{r_1+\cdots +r_t=j-i-(2t-1)\\r_i\ge0}}
\prod_{l=1}^t\dim_k\tilde{H}_{r_l}(\Delta(K_{n_l}^{a_l})_{V_l};k)
\]
obtained from Theorem 3.13 simplifies, via the formula in Theorem 3.14, to 
\[
\sum_{\substack{|V|=j\\V=V_1\sqcup\cdots\sqcup V_t}}
\prod_{l=1}^t{{j_l-1}\choose{a_l-1}}.
\]
This may in turn be written as
\[
\sum_{\substack{j_l\ge0\\j=j_1+\cdots+j_t}}\prod_{l=1}^t
{n_l\choose j_l}{{j_l-1}\choose{a_l-1}}.
\]
Now, if some $j_l\le a_l-1$ the corresponding term is zero. So, we may write 
$j_l=r_l+a_l-1$ where $r_l\ge1$ for $l=1,...,t$. The above expression then becomes
\[
\sum_{\substack{r_l\ge1\\r_1+\cdots+r_t=j-d+t}}\prod_{l=1}^t
{n_l\choose{r_l+a_l-1}}{{r_l+a_l-2}\choose{a_l-1}}.
\] 
Since we know that the resolution is linear, we have the equation $j=i+(d-1)$. Using this in 
the last display we get the formula in the theorem.

\end{proof}

\begin{cor}  The $\mathbb{N}$-graded Betti numbers of 
the $d(a,b)$-complete bipartite hypergraph $K_{n,m}^{d(a,b)}$ may be written as
\[
\beta_{i,j}(K_{n,m}^{d(a,b)})=\sum_{\substack{r+s=i+1\\
r,s\ge1}}{n\choose {r+a-1}}{{r+a-2}\choose {a-1}}{m\choose {s+b-1}}{{s+b-2}\choose 
{b-1}}.
\]
\end{cor}

Furthermore, note that by putting $a=b=1$, we get
\[
\beta_{i,j}(K_{n,m}^{d(1,1)})=
\sum_{\substack{p+q=j\\p,q\ge1}}{n\choose p}{m\choose q}.
\]

Now $K_{n,m}^{d(1,1)}=K_{n,m}$, so we have given another proof of Jacques' formula 
for $\beta_{i,j}(K_{n,m})$.

\begin{cor} Given $K_{\bf n}^{d(\bf a)}$ we have
\[
\beta_i(K_{\bf n}^{d(\bf a)})=
\sum_{\substack{r_1+\cdots +r_t=i+t-1\\
r_i\ge1}}[\prod_{l=1}^t \beta_{r_l, r_l+a_l-1}(K_{n_l}^{a_l})]
\]
\[
\mathrm{pd}(K_{\bf n}^{d(\bf a)})=N-(d-1)
\]
where $j=i+(d-1)$.
\end{cor}

\begin{proof}
The first assertion is clear. If we put $i=N-(d-1)$ in the formula we get
\[
\beta_{N-(d-1)}(K_{\bf n}^{d({\bf a})})=
\prod_{l=1}^t \beta_{n_l-(a_l-1)}(K_{n_l}^{a_l})
\]
which is non zero. At the same time we see that if $i>N-(d-1)$ every term in the sum 
is zero because some factor in every term is zero. 
\end{proof}

{\bf Example:}\, Consider $\mc{H}=K_{3,3,3}^{5(1,1,3)}$. If we denote the set of vertices 
of this hypergraph by $\{a,b,c\}\sqcup \{A,B,C\}\sqcup \{d,e,f\}$,  we get
\[
\mc{E(H)}=\{aAdef,aBdef,aCdef,bAdef,bBdef,bCdef,
cAdef,cBdef,cCdef\}.
\]
The Betti numbers are $\beta_0(\mc{H})=1,\, 
\beta_1(\mc{H})=9,\, \beta_2(\mc{H})=18,\, 
\beta_3(\mc{H})=15,\, \beta_4(\mc{H})=6,\, 
\beta_5(\mc{H})=1$.

\begin{cor} The ring $R/I_{\Delta(K_{\bf n}^{d(\bf a)})^\ast}$ is 
Cohen-Macaulay and we have 
\begin{eqnarray*} \dim\Delta(K_{\bf n}^{d({\bf a})})^\ast&=&N-d-1\\
\dim_R(R/I_{\Delta(K_{\bf n}^{d({\bf a})})^\ast}) &=&N-d\\
\mathrm{pd}(R/I_{\Delta(K_{\bf n}^{d({\bf a})})^\ast})&=&d.
\end{eqnarray*}
\end{cor}

\begin{proof}
The Cohen-Macaulayness follows, for example, from the theorem and the 
Eagon-Reiner theorem. By considering the description of the Alexander dual given in 
Lemma 3.15 the first equation is clear and imply the second. The third is a consequence of 
the fact that $R/I_{\Delta(K_{\bf n}^{d({\bf a})})^\ast}$ is Cohen-Macaulay.

\end{proof}

\begin{prop} The Stanley-Reisner ring $R/I(K_{\bf n}^{d({\bf a})})$ of the
$d(a_1,...,a_t)$-complete multipartite hypergraph on vertex set 
$[n_1]\sqcup\cdots\sqcup [n_t]$ is Cohen-Macaulay precisely when 
$a_s=n_s$ for all $s\in\{1,...,t\}$ but possibly one. This single $a_i$ is such that 
it maximizes the expression $a_i+\sum_{j\neq i,\,j=1}^t n_j$.
\end{prop}

\begin{proof}
Let $I_s\sse [n_s]$. It is necessary and sufficient that at least one set 
$I_i$ satisfy $|I_i|<a_i$, for $I_1\sqcup\cdots\sqcup I_t$ to be a face of 
$\Delta(K_{\bf n}^{d({\bf a})})$. Thus, the dimension of $\Delta(K_{\bf n}^{d({\bf a})})$
is
\[
\max\{a_i-2+\sum_{j\neq i,\,j=1}^t n_j;\, i=1,...,t\}
\]
so
\[
\dim_R(R/I(K_{\bf n}^{d({\bf a})}))=\max\{a_i-1+\sum_{j\neq i,\,j=1}^t n_j\}.
\]
We know that $\t{pd}(K_{\bf n}^{d({\bf a})})=N-(d-1)$, so 
$\t{depth}(R/I(K_{\bf n}^{d({\bf a})}))=d-1$. Now, since by construction 
$d=\sum_{s=1}^t a_s$ we are done.
\end{proof}
Note that this again, in a sense, collapses to an ordinary $d$-complete hypergraph.\\
   
One special, and rather intuitive, way of generalizing the complete bipartite graph, 
is to consider the $d(1,...,1)$-complete mulitpartite hypergraph 
$K_{n_1,...,n_t}^{d(1,...,1)}$. According to the above its $ij$'th Betti number is 
given by
\[
\beta_{i,j}(K_{n_1,...,n_t}^{d(1,...,1)})=
\sum_{\substack{r_1+\cdots +r_t=i+t-1\\r_s\ge1}}
\prod_{l=1}^t{{n_l}\choose {r_l}}.
\]
In \cite{Be}, Berge defines what he calls the $d$-partite complete hypergraph. In our 
language this is just $K_{n_1,...,n_d}^{d(1,...,1)}$, so his definition is a 
special case of ours.

\subsection{The $d(I_1,...,I_t)$-complete hypergraph}

In this section we define another class of complete hypergraphs that actually contains 
all of the previously defined classes of complete hypergraphs. We then show that the 
hypergraphs in this new class have linear resolutions. In this way, one may think that some 
of our previous results are superfluous. We argue that they are not. This is because the main 
part of the results so far is about calculating the Betti numbers. The fact that we have had 
linear resolutions have mainly been used as a computational aid.

In the case of the $d$-complete hypergraph we considered all possible $d$-edges and in the 
case of the $d(a_1,...,a_t)$-complete hypergraph, we considered those in which precisely 
$a_s$ elements came from the vertex set $[n_s]$. We are going to keep the vertex set 
$[n_1]\sqcup\cdots\sqcup [n_t]$ of $N=\sum_{s=1}^t n_s$ vertices, but define another edge set. 
For each $s=1,...,t$ let $I_s$ be an interval $[\alpha_s,\beta_s]$ in $\{0,...,n_s\}$. 
We define the $d(I_1,...,I_t)$-{\bf complete multipartite} hypergraph 
to be the $d$-uniform hypergraph on vertex set 
$[n_1]\sqcup\cdots\sqcup [n_t]$, and with edge set consisting of all $d$-edges 
$I_1(a_1)\sqcup\cdots\sqcup I_t(a_t)$. Here $I_s(a_s)$ is a subset of $[n_s]$ of cardinality 
$a_s\in I_s$ and $d=\sum_{s=1}^t a_s$. We immediately see why 
this generalizes previously considered hypergraphs. If $I_s=\{0,...,n_s\}$ for all 
$s=1,...,t$ we have the $d$-complete hypergraph $K_N^d$. If $I_s=\{0,...,\min\{n_s,d-1\}\}$ 
we obtain the $d$-complete multipartite hypergraph $K_{n_1,...,n_s}^d$. By letting $I_s$ 
consist of only one non zero element for all $s$, we obtain the $d(a_1,...,a_t)$-complete 
hypergraph $K_{n_1,...,n_t}^{d(a_1,...,a_t)}$. 
So, we already know that some special instances of $K_{n_1,...,n_t}^{d(I_1,...,I_t)}$ 
have linear resolutions.\\

One easily realizes that two different sets of intervals $I_1,...,I_t$ and $J_1,...,J_t$ 
say, may yield the same hypergraph. Just consider the case where $I_s=\{a_s\}$ for 
all $s$, $\sum_{s=1}^t a_s=d$, and $J_s=\{a_s\}$ for all $s\neq1$, $J_1=[a_1, a_1+1]$. 
However, to obtain a different hypergraph, we just need to change $J_2$ say, to 
$[a_2-1, a_2]$. 

From now on, without loss of generality, we will assume that the sequence of intervals 
$I_1,...,I_t$ in a hypergraph $K_{n_1,...,n_t}^{d(I_1,...,I_t)}$ satisfies the 
following property: If $I_s=[\alpha_s,\beta_s]$ for $s=1,...,t$, then
\begin{eqnarray*}
\alpha_s+\sum_{j\neq s}\beta_j\ge d\\
\beta_s+\sum_{j\neq s} \alpha_j\le d
\end{eqnarray*}
holds for every $s$. In other words, we assume that there is no redundancy in the sense that 
every element $a_s\in I_s$ is part of an edge in the hypergraph.

It is clear that a set of intervals $I_1,...,I_t$ corresponding to a hypergraph
$K_{n_1,...,n_t}^{d(I_1,...,I_t)}$ can be constructed 
from at least one sequence $a_1,...,a_t$, $d=\sum_{s=1}^t a_s$, corresponding to a 
$d(a_1,...,a_t)$-complete hypergraph, by successively changing the intervals by 
extending one of them (or possibly two of them depending on the situation) in such a way 
that the inequalities above remains true in each step. The following example will clarify 
this idea.\\
\\
{\bf Example:}\, Suppose $a_1+a_2+a_3=d$ and consider $K_{n_1,n_2,n_3}^{d(I_1,I_2,I_3)}$ with 
$I_1=[a_1-1,a_1],\,I_2=[a_2-1,a_2+1],\,I_3=[a_3,a_3+1]$. These intervals can be constructed 
from the trivial intervals $I_1=\{a_1\}$, $I_2=\{a_2\}$ and $I_3=\{a_3\}$ in the 
following way:
\[
\{a_1\},\,\{a_2\},\,\{a_3\}\rightsquigarrow
[a_1-1,a_1],\,[a_2,a_2+1],\,\{a_3\}\rightsquigarrow
[a_1-1,a_1],\,[a_2-1,a_2+1],\,[a_3,a_3+1].
\]

\begin{thm} The Stanley-Reisner ring $R/I_{\Delta(K_{n_1,...,n_t}^{d(I_1,...I_t)})}$ of the 
$d(I_1,...,I_t)$-complete multipartite hypergraph has linear resolution, and 
$\beta_{i,j}(K_{n_1,...,n_t}^{d(I_1,...,I_t)})\neq0$ only if $j=i+(d-1)$.
\end{thm}
\begin{proof}
We start by noting that if $t=1$, then $\Delta(K_{n_1,...,n_t}^{d(I_1,...,I_t)})=
K_{n}^d$ which we know has linear resolution, and non zero homology only in degree $d-2$. 
If $n<d$ we have the $d$-uniform hypergraph with empty edge set and this also has linear 
resolution. Since the set of intervals $I_1,...,I_t$ that corresponds to a hypergraph 
$K_{n_1,...,n_t}^{d(I_1,...,I_t)}$ can be constructed (as above) from several sequences 
$a_1,...,a_t$, $a_1+\cdots +a_t=d$ we may, without loss of generality, assume that 
$I_t=[a_t, a_t+r]$ for some positive integer $r$. Having already gone through the case 
where $t=1$, let us assume that all hypergraphs $K_{n_1,...,n_{t-1}}^{d(I_1,...,I_{t-1})}$ 
have linear resolutions and non zero homology only in degree $d-2$.\\

Given $K_{n_1,...,n_t}^{d(I_1,...,I_t)}$ as above, the expression 
$K_{n_1,...,n_{t-1}}^{(d-a_t-s)(I_1,...,I_{t-1})}$, $s\ge0$ makes sense. 
It means precisely what it says but there is one little problem, the intervals 
$I_1,...,I_{t-1}$ may no longer be as small as possible. There may very well be an element 
$a_j\in I_j$ for some $j=1,...,t-1$, that can not be used in a partition of $d-a_t-s$. 
So, we really should use some other symbols $I'_1,...,I'_{t-1}$ for the intervals, as they may 
depend on $s$. We will however, for convenience, allow this abuse of notation 
in this proof.\\

The next observation we make is that 
\[
\Delta(K_{n_1,...,n_{t-1}}^{(d-a_t-s)(I_1,...,I_{t-1})})\sse
\Delta(K_{n_1,...,n_{t-1}}^{(d-a_t-s+1)(I_1,...,I_{t-1})}).
\]   
Indeed, a face in the first complex can not contain an edge from 
$K_{n_1,...,n_{t-1}}^{(d-a_t-s+1)(I_1,...,I_{t-1})}$, since it would then automatically 
also contain an edge from $K_{n_1,...,n_{t-1}}^{(d-a_t-s)(I_1,...,I_{t-1})}$.\\
By considering the minimal non faces in the complex, one realizes that 
$\Delta(K_{n_1,...,n_t}^{d(I_1,...,I_t)})$ has the following expression.
\[
\bigl[\,\underbrace{\Delta(K_{n_1,...,n_{t-1}}^{(d-a_t)(I_1,...,I_{t-1})})\ast\Delta_{n_t}}_{L_0}
\cup\underbrace{\Delta_{n_1+\cdots +n_{t-1}}\ast\Gamma_{a_t-2}(n_t)}_{M_0}\,\bigr]\,\bigcap
\] 
\[\vdots
\]
\[
\bigl[\,\underbrace{\Delta(K_{n_1,...,n_{t-1}}^{(d-a_t-r+1)(I_1,...,I_{t-1})})\ast
\Delta_{n_t}}_{L_{r-1}}\cup\underbrace{\Delta_{n_1+\cdots +n_{t-1}}\ast
\Gamma_{a_t+r-3}(n_t)}_{M_{r-1}}\,\bigr]\,\bigcap
\] 
\[
\bigl[\,\underbrace{\Delta(K_{n_1,...,n_{t-1}}^{(d-a_t-r)(I_1,...,I_{t-1})})\ast
\Delta_{n_t}}_{L_r}\cup\underbrace{\Delta_{n_1+\cdots +n_{t-1}}\ast
\Gamma_{a_t+r-2}(n_t)}_{M_r}\,\bigr].
\] 
Here the expression in each row correspond to the faces that do not contain one type of 
minimal non face. We take the intersection of these to obtain the faces that do not contain 
any minimal non face, in other words the whole complex.\\
Put $A_0=L_0,\, B_0=M_0,\, \Delta(0)=A_0\cup B_0$ and define recursively 
$A_r=L_r\cap\Delta(r-1),\, B_r=M_r\cap\Delta(r-1),\, \Delta(r)=A_r\cup B_r$. We will now deduce 
explicit formulas for $A_r,\,B_r$ and $A_r\cap B_r$. Having done this, easy use of Mayer-Vietoris 
together with the induction hypothesis will give our result.\\

We start by considering $A_r$.
\[
A_r=L_r\cap\Delta(r-1)=L_r\cap\bigl(A_{r-1}\cup B_{r-1}\bigr)=
L_r\cap\bigl((L_{r-1}\cup M_{r-1})\cap\Delta(r-2)\bigr)=
\]  
\[
L_r\cap\Delta(r-2)=\cdots=L_r\cap\Delta(0)=
\Delta(K_{n_1,...,n_{t-1}}^{(d-a_t-r)(I_1,...,I_{t-1})})\ast\Delta_{n_t}.
\] 
The expression for $B_r$ we will prove by induction.\\
\\
{\bf Claim:}\, $B_r=\bigl[\,\bigcup_{s=0}^{r-1}
\Delta(K_{n_1,...,n_{t-1}}^{(d-a_t-s)(I_1,...,I_{t-1})})\ast\Gamma_{a_t+s-1}(n_t)\,\bigr]
\cup\Delta_{n_1+\cdots+n_{t-1}}\ast\Gamma_{a_t-2}(n_t)$.\\
\\
If we interpret the expression in brackets as $\emptyset$ when $r=0$, it is clear that the 
formula holds for $r=0$. Now 
\[
B_1=\bigl(\Delta_{n_1+\cdots+n_{t-1}}\ast\Gamma_{a_t-1}(n_t)\bigr)
\bigcap\bigl(A_0\cup B_0\bigr)=
\]
\[
\Delta(K_{n_1,...,n_{t-1}}^{(d-a_t)(I_1,...,I_{t-1})})\ast
\Gamma_{a_t-1}(n_t)\cup\Delta_{n_1+\cdots+n_{t-1}}\ast\Gamma_{a_t-2}(n_t)
\]
so the formula holds for 1 as well. 
Assume that the formula holds for $r$. Then
\[
B_{r+1}=M_{r+1}\bigcap\bigl(A_r\cup B_r\bigr).
\]
Let us investigate $M_{r+1}\cap A_r$ and $M_{r+1}\cap B_r$ separately. Using the expression for 
$A_r$ that we already have, we immediately get
\[
M_{r+1}\cap A_r=\Delta(K_{n_1,...,n_{t-1}}^{(d-a_t-r)(I_1,...,I_{t-1})})\ast
\Gamma_{a_t+r-1}(n_t).
\]
Furthermore
\[
M_{r+1}\cap B_r= M_{r+1}\bigcap\bigl[\bigl(
\bigcup_{s=0}^{r-1}\Delta(K_{n_1,...,n_{t-1}}^{(d-a_t-s)(I_1,...,I_{t-1})})\ast
\Gamma_{a_t+s-1}(n_t)\bigr)\bigr]
\cup\Delta_{n_1+\cdots+n_{t-1}}\ast\Gamma_{a_t-2}(n_t)\bigr]=
\]
\[
[\bigl(
\bigcup_{s=0}^{r-1}\Delta(K_{n_1,...,n_{t-1}}^{(d-a_t-s)(I_1,...,I_{t-1})})\ast
\Gamma_{a_t+s-1}(n_t)\bigr)]
\cup\Delta_{n_1+\cdots+n_{t-1}}\ast\Gamma_{a_t-2}(n_t).
\]
So,
\[
B_{r+1}=\bigl(M_{r+1}\cap A_r\bigr)\cup\bigl(M_{r+1}\cap B_r\bigr)=
\]
\[
\bigl[\,\bigcup_{s=0}^{r}
\Delta(K_{n_1,...,n_{t-1}}^{(d-a_t-s)(I_1,...,I_{t-1})})\ast\Gamma_{a_t+s-1}(n_t)\,\bigr]
\cup\Delta_{n_1+\cdots+n_{t-1}}\ast\Gamma_{a_t-2}(n_t)
\]
and we have proved the claim.\\
Lastly, we consider $A_r\cap B_r$. Using the expressions that we have just deduced, we get
\[
A_r\cap B_r=\bigl[\,\Delta(K_{n_1,...,n_{t-1}}^{(d-a_t-r)(I_1,...,I_{t-1})})\ast\Delta_{n_t}
\,\bigr]\bigcap
\]
\[
\bigl[\,\big(\bigcup_{s=0}^{r-1}
\Delta(K_{n_1,...,n_{t-1}}^{(d-a_t-s)(I_1,...,I_{t-1})})\ast\Gamma_{a_t+s-1}(n_t)\,\bigr)
\cup\Delta_{n_1+\cdots+n_{t-1}}\ast\Gamma_{a_t-2}(n_t)\,\bigr]=
\]
\[
\bigl(\,\Delta(K_{n_1,...,n_{t-1}}^{(d-a_t-r)(I_1,...,I_{t-1})})\ast
\Gamma_{a_t+r-2}(n_t)\,\bigr)\bigcup
\bigl(\,\Delta(K_{n_1,...,n_{t-1}}^{(d-a_t-r)(I_1,...,I_{t-1})})\ast
\Gamma_{a_t-2}(n_t)\,\bigr)=
\]
\[
\Delta(K_{n_1,...,n_{t-1}}^{(d-a_t-r)(I_1,...,I_{t-1})})\ast
\Gamma_{a_t+r-2}(n_t).
\]
Now, since we have descriptions of $A_r,\,B_r$ and $A_r\cap B_r$, it will be rather 
easy to finish the proof.\\
$A_r$ is a cone, and hence have no homology at all. For $B_r$ we write
\[
B_r=\underbrace{\bigl[\,\bigcup_{s=0}^{r-1}
\Delta(K_{n_1,...,n_{t-1}}^{(d-a_t-s)(I_1,...,I_{t-1})})\ast\Gamma_{a_t+s-1}(n_t)\,\bigr]}_{L}
\cup\underbrace{\Delta_{n_1+\cdots+n_{t-1}}\ast\Gamma_{a_t-2}(n_t)}_{M}.
\]
$M$ is a cone and have no homology at all, and $L$ can only have homology in degree $d-2$. 
This follows easily by induction on $r$ using Mayer-Vietoris. The case where $r=1$ follows 
directly from the results in section 2.5. Furthermore one easily sees that 
$L\cap M=\Delta(K_{n_1,...,n_{t-1}}^{(d-a_t)(I_1,...,I_{t-1})})\ast\Gamma_{a_t-2}(n_t)$. 
It follows again from the result in section 2.5 that this complex can only have homology in 
degree $d-3$. Thus, using Mayer-Vietoris on $B_r$ we conclude that the homology 
of $B_r$ can be non zero only in degree $d-2$.\\

Using K\"unneth's formula again, we immediately conclude that 
$A_r\cap B_r$ only have homology in degree $d-3$. Thus, the exact sequence 
\[
0\to\mc{C}.(A_r\cap B_r)\to\mc{C}.(A_r)\oplus\mc{C}.(B_r)\to\mc{C}.(A_r\cup B_r)\to0
\]
gives, via Mayer-Vietoris, that the homology of $A_r\cup B_r$ can be non zero only in degree 
$d-2$. Now we are almost done. We aim to use Hochster's formula to conclude that the resolution 
of $K_{n_1,...,n_t}^{d(I_1,...,I_t)}$ is linear. Hence, we would like to conclude 
that the homology of every restriction $\Delta(K_{n_1,...,n_t}^{d(I_1,...,I_t)})_V$ behaves 
precisely like that of the whole complex. In other words, that it only can exist in degree 
$d-2$. If $V\cap [n_s]\neq\emptyset$ for all $s=1,...,t$ and the induced hypergraph
$(K_{n_1,...,n_t}^{d(I_1,...,I_t)})_V$ has non empty edge set, then the restriction 
$\Delta(K_{n_1,...,n_t}^{d(I_1,...,I_t)})_V$ has 
precisley the same form as the original complex $\Delta(K_{n_1,...,n_t}^{d(I_1,...,I_t)})$. 
Thus, by the above, it may only have homology in degree $d-2$. Next assume that 
$V\cap [n_s]=\emptyset$ for at least one $s$. Then 
$F\in\Delta(K_{n_1,...,n_t}^{d(I_1,...,I_t)})_V$ if and only if $\Delta_{n_s}\ast F
\in\Delta(K_{n_1,...,n_t}^{d(I_1,...,I_t)})_V$. Hence, in this case 
$\Delta(K_{n_1,...,n_t}^{d(I_1,...,I_t)})_V$ is a cone and have no homology. 
The last case to consider is if $V\cap [n_s]\neq\emptyset$ for all $s$, but the edge set of 
$(K_{n_1,...,n_t}^{d(I_1,...,I_t)})_V$ is empty. Also in this case we have a cone, in fact, 
a simplex. This is easy since if $(K_{n_1,...,n_t}^{d(I_1,...,I_t)})_V$ has no edges, then 
there are no minimal non faces in $\Delta(K_{n_1,...,n_t}^{d(I_1,...,I_t)})_V$. We have thus 
proved that if
$\dim_k\tilde{H}_{l}(\Delta(K_{n_1,...,n_t}^{d(I_1,...,I_t)})_V;k)\neq0$, then $l=d-2$. Hence 
Hochster's formula gives that $\beta_{i,j}(K_{n_1,...,n_t}^{d(I_1,...,I_t)})\neq0$ only 
if $j=i+(d-1)$ and we are done.

\end{proof}

{\bf Example:}\, Consider $\mc{H}=K_{3,3,3}^{5(I_1,I_2,I_3)}$ where 
$I_1=[1,2], I_2=\{1\}, I_3=[2,3]$. There are 36 5-edges in this hypergraph and 
using a computer one easily computes the Betti numbers. They are
$\beta_0(\mc{H})=1,\, \beta_1(\mc{H})=36,\, \beta_2(\mc{H})=90,\, 
\beta_3(\mc{H})=87,\, \beta_4(\mc{H})=39,\, \beta_5(\mc{H})=7$.\\

By considering the edges in $K_{n_1,...,n_t}^{d(I_1,...,I_t)}$ and the description of the 
Alexander dual of $\Delta(K_{\bf n}^{d(\bf a)})$, we obtain the following description of 
$\Delta(K_{n_1,...,n_t}^{d(I_1,...,I_t)})^\ast$.
\[
\Delta(K_{n_1,...,n_t}^{d(I_1,...,I_t)})^\ast=\bigcup_{\substack{a_1+\cdots +a_t=d\\a_s\in I_s}}
\Gamma_{n_1-a_1-1}(n_1)\ast\cdots\ast\Gamma_{n_t-a_t-1}(n_t).
\]
We immediately get the following
\begin{cor} The ring $R/I_{\Delta(K_{n_1,...,n_t}^{d(I_1,...,I_t)})^\ast}$ is Cohen-Macualay and 
we have
\begin{eqnarray*} \dim\Delta(K_{n_1,...,n_t}^{d(I_1,...,I_t)})^\ast&=&N-d-1\\
\dim_R(R/I_{\Delta(K_{n_1,...,n_t}^{d(I_1,...,I_t)})^\ast})&=&N-d\\
\mathrm{pd}(R/I_{\Delta(K_{n_1,...,n_t}^{d(I_1,...,I_t)})^\ast})&=&d.
\end{eqnarray*}
\end{cor}

\bibliographystyle{plain}
\bibliography{ref}

\end{document}